\pdfoutput=1
\RequirePackage{ifpdf}
\ifpdf 
\documentclass[pdftex]{sigma}
\else
\documentclass{sigma}
\fi

\begin{document}

\newcommand{\arXivNumber}{1205.2946}

\allowdisplaybreaks

\renewcommand{\thefootnote}{$\star$}

\renewcommand{\PaperNumber}{007}

\FirstPageHeading

\ShortArticleName{On a~Certain Subalgebra of $U_q(\widehat{\mathfrak{sl}}_2)$ Related to the Degenerate $q$-Onsager Algebra}

\ArticleName{On a~Certain Subalgebra of $\boldsymbol{U_q(\widehat{\mathfrak{sl}}_2)}$\\
Related to the Degenerate $\boldsymbol{q}$-Onsager Algebra\footnote{This paper is a~contribution to the Special Issue on
New Directions in Lie Theory.
The full collection is available at
\href{http://www.emis.de/journals/SIGMA/LieTheory2014.html}{http://www.emis.de/journals/SIGMA/LieTheory2014.html}}}

\Author{Tomoya HATTAI~$^\dag$ and Tatsuro ITO~$^\ddag$}

\AuthorNameForHeading{T.~Hattai and T.~Ito}

\Address{$^\dag$~Iida Highschool, 1-1, Nonoe, Suzu, Ishikawa 927-1213, Japan}
\EmailD{\href{mailto:tmyhtti@m2.ishikawa-c.ed.jp}{tmyhtti@m2.ishikawa-c.ed.jp}}

\Address{$^\ddag$~School of Mathematical Sciences, Anhui University, 111 Jiulong Road, Hefei 230601, China}
\EmailD{\href{mailto:tito@staff.kanazawa-u.ac.jp}{tito@staff.kanazawa-u.ac.jp}}

\ArticleDates{Received September 30, 2014, in f\/inal form January 15, 2015; Published online January 19, 2015}

\Abstract{In~[\textit{Kyushu~J.~Math.} \textbf{64} (2010), 81--144], it is discussed that a~certain subalgebra of the quantum af\/f\/ine algebra
$U_q(\widehat{\mathfrak{sl}}_2)$ controls the second kind TD-algebra of type~I (the degenerate $q$-Onsager algebra).
The subalgebra, which we denote by~$U'_q(\widehat{\mathfrak{sl}}_2)$, is generated by $e_0^+$, $e_1^\pm$, $k_i^{\pm1}$
$(i=0,1)$ with $e^-_0$ missing from the Chevalley generators~$e_i^\pm$, $k_i^{\pm1}$ $(i=0,1)$
of~$U_q(\widehat{\mathfrak{sl}}_2)$.
In this paper, we determine the f\/inite-dimensional irreducible representations of~$U'_q(\widehat{\mathfrak{sl}}_2)$.
Intertwiners are also determined.}


\Classification{17B37; 05E30}

\renewcommand{\thefootnote}{\arabic{footnote}} 
\setcounter{footnote}{0}

\section{Introduction}

Throughout this paper, the ground f\/ield is $\mathbb{C}$ and~$q$ stands for a~nonzero scalar that is not a~root of unity.
The symbols $\varepsilon$, $\varepsilon^*$ stand for an integer chosen from $\{0,1\}$.
Let $\mathcal{A}_q=\mathcal{A}_q^{(\varepsilon,\varepsilon^*)}$ denote the associative algebra with $1$ generated
by~$z$, $z^*$ subject to the def\/ining relations~\cite{augTD}
\begin{gather*}
\text{(TD)}\quad \begin{cases}
  \big[z,[z,[z,z^*]_q]_{q^{-1}}\big]=-\varepsilon\big(q^2-q^{-2}\big)^2[z,z^*],
\\
  \big[z^*,[z^*,[z^*,z]_q]_{q^{-1}}\big]=-\varepsilon^*\big(q^2-q^{-2}\big)^2[z^*,z],
\end{cases}
\end{gather*}
where $[X,Y]=XY-YX$, $[X,Y]_q=qXY-q^{-1}YX$.
This paper deals with a~subalgebra of the quantum af\/f\/ine algebra $U_q(\widehat{\mathfrak{sl}}_2)$ that is closely
related to $\mathcal{A}_q$ in the case of $(\varepsilon, \varepsilon^*)=(1,0)$.
If $(\varepsilon, \varepsilon^*)=(0,0)$, $\mathcal{A}_q$ is isomorphic to the positive part of
$U_q(\widehat{\mathfrak{sl}}_2)$.
If $(\varepsilon, \varepsilon^*)=(1,1)$, $\mathcal{A}_q$ is called the $q$-Onsager algebra.
If $(\varepsilon, \varepsilon^*)=(1,0)$, $\mathcal{A}_q$ may well be called the {\it degenerate $q$-Onsager algebra}.

\looseness=1
The algebra $\mathcal{A}_q$ arises in the course of the classif\/ication of TD-pairs of type~I, which is a~critically
important step in the study of representations of Terwilliger algebras for $P$-~and $Q$-polynomial
association schemes~\cite{ITT}.
For this reason, $\mathcal{A}_q$ is called the TD-algebra of type~I.
Precisely speaking, the TD-algebra of type~I is
standardized to be the algebra $\mathcal{A}_q$, where~$q$ is the main parameter for TD-pairs of type~I;
so $q^2 \neq \pm 1$ and~$q$ is allowed to be a~root of unity.
In our case where we assume~$q$ is not a~root of unity, the classif\/ication of the TD-pairs of type~I is equivalent to
determining the f\/inite-dimensional irreducible representations $\rho: \mathcal{A}_q\to \textrm{End}(V)$
with the property that $\rho(z)$, $\rho(z^*)$ are both diagonalizable.
Such irreducible representations of~$\mathcal{A}_q$ are determined in~\cite{augTD} via embeddings of $\mathcal{A}_q$
into the augmented TD-algebra $\mathcal{T}_q$.
(In the case of $(\varepsilon,\varepsilon^*)=(1,1)$, the diagonalizability condition of $\rho(z)$, $\rho(z^*)$ can be
dropped, because it turns out that this condition always holds for every f\/inite-dimensional irreducible
representation~$\rho$ of the $q$-Onsager algebra~$\mathcal{A}_q$.) $\mathcal{T}_q$ is easier than~$\mathcal{A}_q$ to
study representations for, and each f\/inite-dimensional irreducible representation $\rho:\mathcal{A}_q\to
\textrm{End}(V)$ with~$\rho(z)$,~$\rho(z^*)$ diagonalizable can be extended to a~f\/inite-dimensional irreducible
representation of~$\mathcal{T}_q$ via a~certain embedding of~$\mathcal{A}_q$ into~$\mathcal{T}_q$.

The augmented TD-algebra $\mathcal{T}_q=\mathcal{T}_q^{(\varepsilon, \varepsilon^*)}$ is the associative algebra with
$1$ generated by $x$, $y$, $k^{\pm1}$ subject to the def\/ining relations
\begin{gather}
\textrm{(TD)}_0
\quad
\begin{cases} kk^{-1}=k^{-1}k=1,
\\
kxk^{-1}=q^2x,
\\
kyk^{-1}=q^{-2}y,
\end{cases}
\label{Eq2}
\end{gather}
and
\begin{gather}
\textrm{(TD)}_1
\quad
\begin{cases}
\big[x,[x,[x,y]_q]_{q^{-1}}\big]=\delta\big(\varepsilon^*x^2k^2-\varepsilon k^{-2}x^2\big),
\\
\big[y,[y,[y,x]_q]_{q^{-1}}\big]=\delta\big({-}\varepsilon^*k^2y^2+\varepsilon y^2k^{-2}\big),
\end{cases}
\label{Eq3}
\end{gather}
where $\delta=-(q-q^{-1})(q^2-q^{-2})(q^3-q^{-3})q^4$.
The f\/inite-dimensional irreducible representations of $\mathcal{T}_q$ are determined in~\cite{augTD} via embeddings of
$\mathcal{T}_q$ into the $U_q(\mathfrak{sl}_2)$-loop algebra $U_q(L(\mathfrak{sl}_2))$.

Let $e_i^\pm$, $k_i^{\pm1}$ $(i=0,1)$
be the Chevalley generators of $U_q(L(\mathfrak{sl}_2))$.
So the def\/ining relations of $U_q(L(\mathfrak{sl}_2))$ are
\begin{gather}
\begin{split}
&k_0k_1=k_1k_0=1,
\qquad
 k_ik^{-1}_i=k^{-1}_ik_i=1,
\qquad
k_ie^\pm_ik^{-1}_i=q^{\pm2}e^\pm_i,
\\
&k_ie^\pm_jk^{-1}_i=q^{\mp2}e^\pm_j,
\quad
i\neq j,
\qquad
 [e^+_i,e^-_i]=\dfrac{k_i-k^{-1}_i}{q-q^{-1}},
\qquad [e^+_i,e^-_j]=0,
\quad
i\neq j,
\\
&\big[e^\pm_i,[e^\pm_i,[e^\pm_i,e^\pm_j]_q]_{q^{-1}}]=0,
\qquad
i\neq j.
\end{split}
\label{Eq4}
\end{gather}
Note that if $k_0k_1=k_1k_0=1$ is replaced by $k_0k_1=k_1k_0$ in~\eqref{Eq4}, we have the quantum af\/f\/ine algebra
$U_q(\widehat{\mathfrak{sl}}_2)$: $U_q(L(\mathfrak{sl}_2))$ is the quotient algebra of $U_q(\widehat{\mathfrak{sl}}_2)$
by the two-sided ideal generated by $k_0k_1-1$.
For a~nonzero scalar~$s$, def\/ine the elements $x(s)$, $y(s)$, $k(s)$ of $U_q(L(\mathfrak{sl}_2))$~by
\begin{gather}
\begin{split}
&x(s)=-q^{-1}\big(q-q^{-1}\big)^2\big(se_0^++\varepsilon s^{-1}e_1^-k_1\big),
\\
&y(s)=\varepsilon^*se_0^-k_0+s^{-1}e_1^+,
\\
& k(s)=sk_0.
\end{split}\label{Eq5}
\end{gather}
Then the mapping
\begin{gather}
\varphi_s: \ \mathcal{T}_q\to U_q(L(\mathfrak{sl}_2)),
\qquad
x,y,k\mapsto x(s), y(s), k(s),
\label{Eq6}
\end{gather}
gives an injective algebra homomorphism.
If $(\varepsilon,\varepsilon^*)=(0,0)$, the image $\varphi_s(\mathcal{T}_q)$ coincides with the Borel subalgebra
generated by~$e_i^+$,
$k_i^{\pm1}$
$(i=0,1)$.
If $(\varepsilon,\varepsilon^*)=(1,0)$, the image $\varphi_s(\mathcal{T}_q)$ is properly contained in the subalgebra
generated by $e_0^+$, $e_1^\pm$, $k_i^{\pm1}$
$(i=0,1)$ with $e_0^-$ missing from the generators; we denote this subalgebra by $U'_q(L(\mathfrak{sl}_2))$.
Through the natural homomorphism $U_q(\widehat{\mathfrak{sl}}_2)\to U_q(L(\mathfrak{sl}_2))$, pull back the
subalgebra $U'_q(L(\mathfrak{sl}_2))$ and denote the pre-image by $U'_q(\widehat{\mathfrak{sl}}_2)$:
\begin{gather*}
U'_q(\widehat{\mathfrak{sl}}_2)=\big\langle e_0^+,e_1^\pm,k_i^{\pm1}
\,
|
\,
i=0,1\big\rangle \subset U_q(\widehat{\mathfrak{sl}}_2).
\end{gather*}

\looseness=1
In~\cite{augTD}, it is shown that in the case of $(\varepsilon,\varepsilon^*)=(1,0)$, all the f\/inite-dimensional
irreducible representations of $\mathcal{T}_q$ are produced by tensor products of evaluation modules for
$U'_q(L(\mathfrak{sl}_2))$ via the embedding $\varphi_s$ of $\mathcal{T}_q$ into $U'_q(L(\mathfrak{sl}_2))$.
Using this fact and the Drinfel'd polynomials, we show in this paper that there are no other f\/inite-dimensional
irreducible representations of $U'_q(L(\mathfrak{sl}_2))$ and hence of $U'_q(\widehat{\mathfrak{sl}}_2)$ than those
af\/forded by tensor products of evaluation modules, if we apply suitable automorphisms of $U'_q(L(\mathfrak{sl}_2))$,
$U'_q(\widehat{\mathfrak{sl}}_2)$ to adjust the types of the representations to be~$(1, 1)$.
Here we note that the evaluation parameters are allowed to be zero for $U'_q(L(\mathfrak{sl}_2))$,
$U'_q(\widehat{\mathfrak{sl}}_2)$.
Details will be discussed in Sections~\ref{Section2} and~\ref{Section3}, where the isomorphism classes of f\/inite-dimensional irreducible
representations of $U'_q(\widehat{\mathfrak{sl}}_2)$ are also determined.
In Section~\ref{Section4}, intertwiners will be determined for f\/inite-dimensional irreducible
$U'_q(\widehat{\mathfrak{sl}}_2)$-modules.

In our approach, Drinfel'd polynomials are the key tool for the classif\/ication of f\/inite-dimensional irreducible
representations of $U_q(\widehat{\mathfrak{sl}}_2)$, $U'_q(\widehat{\mathfrak{sl}}_2)$, although they are not the main
subject of this paper.
They are def\/ined in~\cite{augTD}, and the point is that they are directly attached to $\mathcal{T}_q$-modules, not to
$U_q(\widehat{\mathfrak{sl}}_2)$- or $U'_q(\widehat{\mathfrak{sl}}_2)$-modules.
(In the case of $(\varepsilon,\varepsilon^*)=(0,0)$, they turn out to coincide with the original ones up to the
reciprocal of the variable.) So in our approach to the case of $(\varepsilon,\varepsilon^*)=(0,0)$, f\/inite-dimensional
irreducible representations are naturally classif\/ied f\/irstly for the Borel subalgebra of
$U_q(\widehat{\mathfrak{sl}}_2)$ and then for $U_q(\widehat{\mathfrak{sl}}_2)$ itself.
This will be brief\/ly demonstrated in Section~\ref{Section3} as a~warm-up for the case of $(\varepsilon,\varepsilon^*)=(1,0)$, thus
giving another proof to the classical classif\/ication theorem of Chari--Pressley~\cite{CP} and to the main theorems
(Theorems~1.16 and~1.17) of~\cite{BT}.

We now review Drinfel'd polynomials for $\mathcal{T}_q$-modules \cite[p.~119]{augTD}.
Let~$V$ be a~f\/inite-dimensio\-nal $\mathcal{T}_q$-module.
We assume the following properties for~$V$:
\begin{enumerate}\itemsep=0pt
\item[]$({\rm D})_0$:  $k$ is diagonalizable on $V$ with $V=\bigoplus\limits_{i=0}^d U_i$,
$k|_{U_i}=s q^{2i-d}$, $0 \leq i \leq d$,  for some nonzero constant $s$;
\item[]$({\rm D})_1$: $\dim U_0=1$.
\end{enumerate}
By the relations $({\rm TD})_0: k k^{-1}=k^{-1}k=1$, $kxk^{-1}=q^2x$, $kyk^{-1}=q^{-2}y$, it holds that $x U_i
\subseteq U_{i+1}$, $y U_i \subseteq U_{i-1}$ $(0 \leq i \leq d)$, where $U_{-1}=U_{d+1}=0$.
So the one-dimensional subspa\-ce~$U_0$ is invariant under $y^ix^i$ and we have the sequence $\{\sigma_i
\}_{i=0}^{\infty}$ of eigenvalues $\sigma_i$ of $y^ix^i$ on $U_0$: $\sigma_i = y^ix^i|_{U_0}$.
Notice that $\sigma_0=1$ and $\sigma_i=0$, $d+1 \leq i$.
The Drinfel'd polynomial~$P_V(\lambda)$ of the $\mathcal{T}_q$-module~$V$ is def\/ined~by
\begin{gather*}
P_V(\lambda) = \sum\limits^d_{i=0}(-1)^i \frac{\sigma_i}{(q-q^{-1})^{2i}([i]!)^2}
\prod\limits^d_{j=i+1}\big(\lambda-\varepsilon s^{-2} q^{2(d-j)}-\varepsilon^* s^2 q^{-2(d-j)}\big),
\end{gather*}
where $[i]=[i]_q=(q^i-q^{-i})/(q-q^{-1})$ and $[i]!=[1][2]\cdots[i]$ with the understanding of $[0]!=1$.
Since $\sigma_0=1$, $P_V(\lambda)$ is a~monic polynomial of degree~$d$.

If~$V$ is an irreducible $\mathcal{T}_q$-module, it is known that~$V$ in fact satisf\/ies the properties $({\rm D})_0$,
$({\rm D})_1$ \cite[Lemma~1.2, Theorem~1.8]{augTD}, and these properties provide a~rather simple short proof for the
`injective' part of~\cite[Theorem~1.9]{augTD}, i.e., for the fact that the isomorphism class of the irreducible
$\mathcal{T}_q$-module~$V$ is determined by the trio $(\{\sigma_i \}_{i=0}^{\infty},s,d)$.

If~$V$ is a~tensor product of evaluation modules for $U_q(L(\mathfrak{sl}_2))$ in the case of $(\varepsilon,\varepsilon^*)=(1,1),(0,0)$
or for $U'_q(L(\mathfrak{sl}_2))$ in the case of $(\varepsilon, \varepsilon^*)=(1,0)$,
we regard~$V$ as a~$\mathcal{T}_q$-module via the embedding $\varphi_s$ of~\eqref{Eq6}.
Then it is apparent that the $\mathcal{T}_q$-module~$V$ satisf\/ies the properties $({\rm D})_0$, $({\rm D})_1$.
Moreover it is known that a~product formula holds for the Drinfel'd polynomial $P_V(\lambda)$ and it turns out that
$P_V(\lambda)$ does not depend on the parameter~$s$ of the embedding $\varphi_s$~\cite[Theorem~5.2]{augTD}.
The `surjective' part of~\cite[Theorem~1.9]{augTD} follows from the structure of the zeros of the Drinfel'd polynomial
for such a~tensor product of evaluation modules regarded as a~$\mathcal{T}_q$-module via the embedding $\varphi_s$.

\section[Finite-dimensional irreducible representations of $U'_q(\widehat{\mathfrak{sl}}_2)$]{Finite-dimensional
irreducible representations of $\boldsymbol{U'_q(\widehat{\mathfrak{sl}}_2)}$}\label{Section2}

The subalgebra $U'_q(\widehat{\mathfrak{sl}}_2)$ of the quantum af\/f\/ine algebra $U_q(\widehat{\mathfrak{sl}}_2)$ is
generated by $e_0^+$, $e_1^\pm$, $k_i^{\pm1}$
$(i=0,1)$, $e_0^-$ missing from the generators, and has by the triangular decomposition of
$U_q(\widehat{\mathfrak{sl}}_2)$ the def\/ining relations
\begin{gather}
\begin{split}
&k_0k_1=k_1k_0,
\qquad
k_ik^{-1}_i=k^{-1}_ik_i=1,
\qquad
k_0e_0^+k_0^{-1}=q^2e_0^+,
\qquad
k_1e^\pm_1k^{-1}_1=q^{\pm2}e^\pm_1,
\\
& k_1e_0^+k_1^{-1}=q^{-2}e_0^+,
\qquad
k_0e^\pm_1k^{-1}_0=q^{\mp2}e^\pm_1,
\qquad
[e^+_1,e^-_1]=\dfrac{k_1-k^{-1}_1}{q-q^{-1}},
\\
& [e^+_0,e^-_1]=0,
\qquad
\big[e^+_i,[e^+_i,[e^+_i,e^+_j]_q]_{q^{-1}}\big]=0,
\quad
i\neq j.
\end{split}
\label{Eq8}
\end{gather}
Note that if $k_0k_1=k_1k_0$ is replaced by $k_0k_1=k_1k_0=1$ in~\eqref{Eq8}, we have the def\/ining relations for
$U'_q(L(\mathfrak{sl}_2))$.

Let~$V$ be a~f\/inite-dimensional irreducible $U'_q(\widehat{\mathfrak{sl}}_2)$-module.
Let us f\/irst observe that the $U'_q(\widehat{\mathfrak{sl}}_2)$-module~$V$ is obtained from
a~$U'_q(L(\mathfrak{sl}_2))$-module by applying an automorphism of $U'_q(\widehat{\mathfrak{sl}}_2)$ as follows.
Since the element $k_0k_1$ belongs to the centre of $U'_q(\widehat{\mathfrak{sl}}_2)$, $k_0k_1$ acts on~$V$ as
a~scalar~$s$ by Schur's lemma.
Since $k_0k_1$ is invertible, the scalar~$s$ is nonzero: $k_0k_1|_V=s\in\mathbb{C}^\times$.
Observe that $U'_q(\widehat{\mathfrak{sl}}_2)$ admits an automorphism that sends $k_0$ to $s^{-1}k_0$ and preserves
$k_1$.
Hence we may assume $k_0k_1|_V=1$.
Then we can regard~$V$ as a~$U'_q(L(\mathfrak{sl}_2))$-module.

Let~$V$ be a~f\/inite-dimensional irreducible $U'_q(L(\mathfrak{sl}_2))$-module.
For a~scalar~$\theta$, set $V(\theta)=\{v\in V\,|\,  k_0v=\theta v\}$.
So if $V(\theta)\neq0$, $\theta$ is an eigenvalue of $k_0$ and $V(\theta)$ is the corresponding eigenspace of~$k_0$.
For an eigenvalue~$\theta$ and an eigenvector $v\in V(\theta)$, it holds that $e_0^+v\in V(q^2\theta)$ by the relation
$k_0e_0^+=q^2e_0^+k_0$ and $e_1^\pm v\in V(q^{\mp2}\theta)$ by $k_0e_1^\pm=q^{\mp2}e_1^\pm k_0$.
We have
\begin{gather}
e_0^+V(\theta)\subseteq V\big(q^2\theta\big),
\qquad
e_1^\pm V(\theta)\subseteq V\big(q^{\mp2}\theta\big).
\label{Eq9}
\end{gather}
If $\dim V=1$, then $e_0^+V=0$, $e_1^\pm V=0$ by~\eqref{Eq9} and $k_0|_V=\pm1$ by
$[e^+_1,e^-_1]=(k_1-k^{-1}_1)/(q-q^{-1})=(k_0^{-1}-k_0)/(q-q^{-1})$.
Such a~$U'_q(L(\mathfrak{sl}_2))$-module~$V$ is said to be \textit{trivial}.
Assume $\dim V\geq 2$.
Choose an eigenvalue~$\theta$ of $k_0$ on~$V$.
Then $\sum\limits_{i\in\mathbb{Z}}V(q^{\pm2i}\theta)$ is invariant under the actions of the generators by~\eqref{Eq9}, and so we
have $V=\sum\limits_{i\in\mathbb{Z}}V(q^{\pm2i}\theta)$ by the irreducibility of the
$U'_q(L(\mathfrak{sl}_2))$-module~$V$.
Since~$V$ is f\/inite-dimensional, there exists a~positive integer~$d$ and a~nonzero scalar~$s_0$ such that the eigenspace
decomposition of~$k_0$ is
\begin{gather}
V=\bigoplus_{i=0}^dV\big(s_0q^{2i-d}\big).
\label{Eq10}
\end{gather}
We want to show that $s_0=\pm1$ holds in~\eqref{Eq10}.

Consider the subalgebra of $U'_q(L(\mathfrak{sl}_2))$ generated by $e_1^\pm$, $k_1^{\pm1}$ and denote it by
$\mathcal{U}:\mathcal{U}=\langle e_1^\pm$, $k_1^{\pm1}\rangle$.
Regard~$V$ as a~$\mathcal{U}$-module.
Since $\mathcal{U}$~is isomorphic to the quantum algebra $U_q(\mathfrak{sl}_2)$, $V$~is a~direct sum of irreducible
$\mathcal{U}$-modules, and for each irreducible $\mathcal{U}$-submodule~$W$ of~$V$, the eigenvalues of $k_1=k_0^{-1}$
on~$W$ are either $\{q^{2i-\ell}\,|\,  0\leq i\leq\ell\}$ or $\{-q^{2i-\ell}\,|\,  0\leq i\leq\ell\}$ for some nonnegative
integer~$\ell$.
This implies that~(i) $s_0=\pm q^m$ for some $m \in \mathbb{Z}$ and~(ii) if~$\theta$ is an eigenvalue of $k_0$, so is
$\theta^{-1}$.
It follows from~(i) that $V=\bigoplus\limits_{i=0}^dV(\pm q^{2i-d+m})$, and so by~(ii), we obtain $m=0$, i.e., $s_0=\pm1$.

Observe that $U'_q(L(\mathfrak{sl}_2))$ admits an automorphism that sends~$k_i$ to~$-k_i$ $(i=0,1)$ and~$e_1^+$ to~$-e_1^+$.
Hence we may assume $s_0=1$ in~\eqref{Eq10}.
Note that in this case, $k_1$ has the eigenvalues $\{s_1q^{2i-d}\,|\,  0\leq i\leq d\}$ with $s_1=1$.
Such an irreducible module or the irreducible representation af\/forded by such is said to be of {\it type} $(1,1)$,
indicating $(s_0,s_1)=(1,1)$.
We conclude that the determination of f\/inite-dimensional irreducible representations for
$U'_q(\widehat{\mathfrak{sl}}_2)$ is, via automorphisms, reduced to that of type $(1,1)$ for $U'_q(L(\mathfrak{sl}_2))$.

In the rest of this section, we shall introduce evaluation modules for $U'_q(L(\mathfrak{sl}_2))$ and state our main
theorem that every f\/inite-dimensional irreducible representation of type $(1,1)$ of $U'_q(L(\mathfrak{sl}_2))$ is
af\/forded by a~tensor product of them.
For $a\in\mathbb{C}$ and $\ell\in\mathbb{Z}_{\geq 0}$, let $V(\ell,a)$ denote the $(\ell+1)$-dimensional vector space
with a~basis $v_0, v_1, \ldots, v_\ell$.
Using~\eqref{Eq8}, it can be routinely verif\/ied that $U'_q(L(\mathfrak{sl}_2))$ acts on $V(\ell,a)$~by
\begin{gather}
\begin{split}
& k_0v_i=q^{2i-\ell}v_i,
\qquad
k_1v_i=q^{\ell-2i}v_i,
\qquad
e^+_0v_i=aq[i+1]v_{i+1},
\\
& e^+_1v_i=[\ell-i+1]v_{i-1},
\qquad
e^-_1v_i=[i+1]v_{i+1},
\end{split}
\label{Eq11}
\end{gather}
where $v_{-1}=v_{\ell+1}=0$ and $[t]=[t]_q=(q^t-q^{-t})/(q-q^{-1})$.
This $U'_q(L(\mathfrak{sl}_2))$-module $V(\ell,a)$ is irreducible and called an \textit{evaluation module}.
The basis $v_0, v_1, \ldots, v_\ell$ of the $U'_q(L(\mathfrak{sl}_2))$-modu\-le~$V(\ell,a)$ is called
a~\textit{standard basis}.
The vector $v_0$ is called the \textit{highest weight vector}.
Note that the evaluation parameter~$a$ is allowed to be zero.
Also note that if~$\ell=0$, $V(\ell,a)$ is the trivial module.
We denote the evaluation module~$V(\ell,0)$ by~$V(\ell)$, allowing $\ell=0$, and use the nota\-tion~$V(\ell,a)$ only for
an evaluation module with $a\neq0$ and $\ell\geq1$.

The $U_q(\mathfrak{sl}_2)$-loop algebra $U_q(L(\mathfrak{sl}_2))$ has the coproduct $\Delta:
U_q(L(\mathfrak{sl}_2))\to U_q(L(\mathfrak{sl}_2))\otimes U_q(L(\mathfrak{sl}_2))$ def\/ined~by
\begin{gather}
\begin{split}
& \Delta(k_i^{\pm1})=k_i^{\pm1}\otimes k_i^{\pm1},
\qquad
\Delta(e_i^+)=k_i\otimes e_i^++e_i^+\otimes 1,
\\
& \Delta(e_i^-k_i)=k_i\otimes e_i^-k_i+e_i^-k_i\otimes 1.
\end{split}
\label{Eq12}
\end{gather}
The subalgebra $U'_q(L(\mathfrak{sl}_2))$ is closed under~$\Delta$.
Thus given a~set of evaluation modules $V(\ell)$, $V(\ell_i,a_i)$ $(1\leq i\leq n)$ for $U'_q(L(\mathfrak{sl}_2))$, the
tensor product
\begin{gather}
V(\ell)\otimes V(\ell_1,a_1)\otimes\dots\otimes V(\ell_n,a_n)
\label{Eq13}
\end{gather}
becomes a~$U'_q(L(\mathfrak{sl}_2))$-module via~$\Delta$.
Note that by the coassociativity of~$\Delta$, the way of putting parentheses in the tensor product of~\eqref{Eq13} does not
af\/fect the isomorphism class as a~$U'_q(L(\mathfrak{sl}_2))$-module.
Also note that if $\ell=0$ in~\eqref{Eq13}, then $V(0)$ is the trivial module and the tensor product of~\eqref{Eq13} is isomorphic to
$V(\ell_1,a_1)\otimes\dots\otimes V(\ell_n,a_n)$ as $U'_q(L(\mathfrak{sl}_2))$-modules.
Finally we allow $n=0$, in which case we understand that the tensor product of~\eqref{Eq13} means $V(\ell)$.

With the evaluation module $V(\ell,a)$, we associate the set $S(\ell,a)$ of scalars
$aq^{-\ell+1}$, $aq^{-\ell+3},\ldots$, $aq^{\ell-1}$:
\begin{gather*}
S(\ell,a)=\big\{aq^{2i-\ell+1}\,|\,  0\leq i\leq\ell-1\big\}.
\end{gather*}
The set $S(\ell,a)$ is called a~\textit{$q$-string} of length~$\ell$.
Two $q$-strings $S(\ell,a)$, $S(\ell',a')$ are said to be \textit{in general position} if either
\begin{enumerate}\itemsep=0pt
\item[(i)] the union $S(\ell,a)\cup S(\ell',a')$ is not a~$q$-string, or
\item[(ii)] one of $S(\ell,a)$, $S(\ell',a')$ includes the other.
\end{enumerate}
Below is the main theorem of this paper.
It classif\/ies the isomorphism classes of the f\/inite-dimensional irreducible $U'_q(L(\mathfrak{sl}_2))$-modules of type $(1,1)$.
\begin{theorem}\label{Theorem1}
The following~$(i)$,~$(ii)$,~$(iii)$,~$(iv)$ hold.
\begin{enumerate}\itemsep=0pt
\item[$(i)$] A~tensor product $V=V(\ell)\otimes V(\ell_1,a_1)\otimes\cdots\otimes V(\ell_n,a_n)$ of evaluation modules is
irreducible as a~$U'_q(L(\mathfrak{sl}_2))$-module if and only if $S(\ell_i,a_i)$, $S(\ell_j,a_j)$ are in general
position for all $i,j\in\{1,2,\ldots,n\}$.
In this case, $V$ is of type $(1,1)$.
\item[$(ii)$] Consider two tensor products $V=V(\ell)\otimes V(\ell_1,a_1)\otimes\cdots\otimes V(\ell_n,a_n)$,
$V'=V(\ell')\otimes V(\ell'_1,a'_1)\otimes\dots\otimes V(\ell'_{m},a'_{m})$ of evaluation modules and assume that they
are both irreducible as a~$U'_q(L(\mathfrak{sl}_2))$-module.
Then $V$, $V'$ are isomorphic as $U'_q(L(\mathfrak{sl}_2))$-modules if and only if $\ell=\ell'$, $n=m$ and
$(\ell_i,a_i)=(\ell'_i,a'_i)$ for all $i$, $1\leq i\leq n$, with a~suitable reordering of the evaluation modules $V(\ell_1,a_1),\ldots,V(\ell_n,a_n)$.
\item[$(iii)$] Every non-trivial finite-dimensional irreducible $U'_q(L(\mathfrak{sl}_2))$-module of type $(1,1)$ is
isomorphic to some tensor product $V(\ell)\otimes V(\ell_1,a_1)\otimes\cdots\otimes V(\ell_n,a_n)$ of evaluation modules.
\item[$(iv)$] If a~tensor product $V=V(\ell)\otimes V(\ell_1,a_1)\otimes\cdots\otimes V(\ell_n,a_n)$ of evaluation modules
is irreducible as a~$U'_q(L(\mathfrak{sl}_2))$-module, then any change of the orderings of the evaluation modules
$V(\ell)$, $V(\ell_1,a_1),\dots, V(\ell_n,a_n)$ for the tensor product does not change the isomorphism class of the
$U'_q(L(\mathfrak{sl}_2))$-module~$V$.
\end{enumerate}
\end{theorem}

\section{Proof of Theorem~\ref{Theorem1}(i),~(ii),~(iii)}\label{Section3}

Discard the evaluation module $V(\ell)$ from the statement of Theorem~\ref{Theorem1} and replace $U'_q(L(\mathfrak{sl}_2))$ by
$U_q(L(\mathfrak{sl}_2))$ or by $\mathcal{B}$, where $\mathcal{B}$ is the Borel subalgebra of $U_q(L(\mathfrak{sl}_2))$
generated by $e_i^+$, $k_i^{\pm1}$
$(i=0,1)$.
Then it precisely describes the classif\/ication of the isomorphism classes of f\/inite-dimensional irreducible modules of
type $(1,1)$ for $U_q(L(\mathfrak{sl}_2))$~\cite{CP} or for $\mathcal{B}$~\cite{BT}.
There is a~one-to-one correspondence of f\/inite-dimensional irreducible modules of type $(1,1)$ between
$U_q(L(\mathfrak{sl}_2))$ and $\mathcal{B}$: every f\/inite-dimensional irreducible $U_q(L(\mathfrak{sl}_2))$-module of
type $(1,1)$ is irreducible as a~$\mathcal{B}$-module and every f\/inite-dimensional irreducible $\mathcal{B}$-module of
type $(1,1)$ is uniquely extended to a~$U_q(L(\mathfrak{sl}_2))$-module.
This sort of correspondence of f\/inite-dimensional irreducible modules exists between $U'_q(L(\mathfrak{sl}_2))$ and
$\mathcal{T}_q$ via the embedding $\varphi_s$ of~\eqref{Eq6}, where $\mathcal{T}_q$ is the augmented TD-algebra with
$(\varepsilon, \varepsilon^*)=(1,0)$.
Namely, we shall show that~(i) every f\/inite-dimensional irreducible $U'_q(L(\mathfrak{sl}_2))$-module of type $(1,1)$ is
irreducible as a~$\mathcal{T}_q$-module via certain embedding $\varphi_s$ of~\eqref{Eq6}, and~(ii) every f\/inite-dimensional
irreducible $\mathcal{T}_q$-module is uniquely extended to a~$U'_q(L(\mathfrak{sl}_2))$-module of type $(1,1)$ via the
embedding $\varphi_s$ of~\eqref{Eq6} with~$s$ uniquely determined.
Since f\/inite-dimensional irreducible $\mathcal{T}_q$-modules are classif\/ied in~\cite{augTD}, this gives a~proof of
Theorem~\ref{Theorem1}.

Apart from the Drinfel'd polynomials, the key to our understanding of the correspondence is the following two lemmas
about $U_q(\mathfrak{sl}_2)$-modules.
Let $\mathcal{U}$ denote the quantum algebra $U_q(\mathfrak{sl}_2)$: $\mathcal{U}$ is the associative algebra with $1$
generated by $X^\pm$, $K^{\pm 1}$ subject to the def\/ining relations
\begin{gather}
 KK^{-1}=K^{-1}K=1,
\qquad
KX^\pm K^{-1}=q^{\pm2}X^\pm,
\qquad
[X^+,X^-]=\dfrac{K-K^{-1}}{q-q^{-1}}.
\label{Eq15}
\end{gather}

\begin{lemma}[\protect{\cite[Lemma 7.5]{augTD}}]\label{Lemma1}
Let~$V$ be a~finite-dimensional $\mathcal{U}$-module that has the following weight-space
$(K$-eigenspace$)$ decomposition:
\begin{gather*}
V=\bigoplus_{i=0}^{d}U_i,
\qquad
K|_{U_i}=q^{2i-d},
\qquad
0\leq i\leq d.
\end{gather*}
Let~$W$ be a~subspace of~$V$ as a~vector space.
Assume that~$W$ is invariant under the actions of~$X^+$ and~$K$:
\begin{gather*}
X^+W\subseteq W,
\qquad
KW\subseteq W.
\end{gather*}
If it holds that
\begin{gather*}
\dim(W\cap U_i)=\dim(W\cap U_{d-i}),
\qquad
0\leq i\leq d,
\end{gather*}
then $X^-W\subseteq W$, i.e., $W$ is a~$\mathcal{U}$-submodule.
\end{lemma}

\begin{lemma}\label{Lemma2}
If~$V$ is a~finite-dimensional~$\mathcal{U}$-module, the action of~$X^-$ on~$V$ is uniquely determined by those of~$X^+$, $K^{\pm1}$ on~$V$.
\end{lemma}

\begin{proof}
The claim holds if~$V$ is irreducible as a~$\mathcal{U}$-module.
By the semi-simplicity of $\mathcal{U}$, it holds generally.
\end{proof}

As a~warm-up for the proof of Theorem~\ref{Theorem1}, we shall demonstrate how to use these lemmas in the case of the corresponding
theorem~\cite{CP} for $U_q(L(\mathfrak{sl}_2))$.
We want, and it is enough, to show part~(iii) of the theorem for $U_q(L(\mathfrak{sl}_2))$ by using the classif\/ication
of f\/inite-dimensional irreducible $\mathcal{B}$-modules.
This is because the parts~(i),~(ii),~(iv) are well-known in advance of~\cite{CP}, while the f\/inite-dimensional
irreducible $\mathcal{B}$-modules are classif\/ied in~\cite{augTD} rather straightforward by the product formula of
Drinfel'd polynomials without using the part~(iii) in question.

Let~$V$ be a~f\/inite-dimensional irreducible $U_q(L(\mathfrak{sl}_2))$-module of type $(1,1)$.
Then $V$ has the weight-space decomposition
\begin{gather*}
V=\bigoplus_{i=0}^{d}U_i,
\qquad
k_0|_{U_i}=q^{2i-d},
\qquad
0\leq i\leq d.
\end{gather*}
Regard~$V$ as a~$\mathcal{B}$-module.
Let~$W$ be a~minimal~$\mathcal{B}$-submodule of~$V$.
Note that~$W$ is irreducible as a~$\mathcal{B}$-module.
We want to show $W=V$, i.e., $V$ is irreducible as a~$\mathcal{B}$-module.
Since the mapping $(e_0^+)^{d-2i}$: $U_i\rightarrow U_{d-i}$ is a~bijection and $W\cap U_i$ is mapped into $W\cap
U_{d-i}$ by $(e_0^+)^{d-2i}$, we have $\dim(W\cap U_i)\leq \dim(W\cap U_{d-i})$,
$0\leq i\leq[d/2]$.
Similarly from the bijection $(e_1^+)^{d-2i}$: $U_{d-i}\rightarrow U_{i}$, we get $\dim(W\cap U_{d-i})\leq\dim(W\cap U_i)$.
Thus it holds that
\begin{gather*}
\dim(W\cap U_i)=\dim(W\cap U_{d-i}),
\qquad
0\leq i\leq d.
\end{gather*}
Consider the algebra homomorphism from $\mathcal{U}$ to $U_q(L(\mathfrak{sl}_2))$ that sends $X^+$, $X^-$, $K^{\pm1}$ to
$e_0^+$, $e_0^-$, $k_0^{\pm 1}$, respectively.
Regard~$V$ as a~$\mathcal{U}$-module via this algebra homomorphism.
Then $X^+W\subseteq W$, $KW\subseteq W$.
Since $\dim(W\cap U_i)=\dim(W\cap U_{d-i})$, $0\leq i\leq d$, we have by Lemma~\ref{Lemma1} that $X^-W\subseteq W$, i.e., $e_0^-W\subseteq W$.
Similarly, Lemma~\ref{Lemma1} can be applied to the $\mathcal{U}$-module~$V$ via the algebra homomorphism from $\mathcal{U}$ to
$U_q(L(\mathfrak{sl}_2))$ that sends $X^+$, $X^-$, $K^{\pm1}$ to $e_1^+$, $e_1^-$, $k_1^{\pm 1}$, respectively, in which
case the weight-space decomposition of the $\mathcal{U}$-module~$V$ is $V=\bigoplus\limits_{i=0}^{d}U_{d-i}$, $
K|_{U_{d-i}}=q^{2i-d}$, $0\leq i\leq d$.
Consequently, we get $X^-W\subseteq W$, i.e., $e_1^-W\subseteq W$.
Thus~$W$ is $U_q(L(\mathfrak{sl}_2))$-invariant and we have $W=V$ by the irreducibility of the
$U_q(L(\mathfrak{sl}_2))$-module~$V$.
We conclude that every f\/inite-dimensional irreducible $U_q(L(\mathfrak{sl}_2))$-module of type $(1,1)$ is irreducible as
a~$\mathcal{B}$-module.

Now consider the class of f\/inite-dimensional irreducible $\mathcal{B}$-modules~$V$, where~$V$ runs through all tensor
products of evaluation modules that are irreducible as a~$U_q(L(\mathfrak{sl}_2))$-module:
\begin{gather*}
V=V(\ell_1,a_1)\otimes\cdots\otimes V(\ell_n,a_n).
\end{gather*}
Then it turns out that the Drinfel'd polynomials $P_V(\lambda)$ of the irreducible $\mathcal{B}$-modules~$V$ exhaust all
that are possible for f\/inite-dimensional irreducible $\mathcal{B}$-modules of type $(1,1)$, as shown in~\cite[Theorem~5.2]{augTD}
by the product formula
\begin{gather*}
P_V(\lambda)=\prod\limits_{i=1}^nP_{V(\ell_i,a_i)}(\lambda),
\qquad
P_{V(\ell_i,a_i)}(\lambda)=\prod\limits_{\zeta\in S(\ell_i,a_i)}(\lambda+\zeta).
\end{gather*}
Since the Drinfel'd polynomial $P_V(\lambda)$ determines the isomorphism class of the $\mathcal{B}$-module~$V$ of type
$(1,1)$ \cite[the injectivity part of Theorem 1.9$'$]{augTD},
there are no other f\/inite-dimensional irreducible
$\mathcal{B}$-modules of type $(1,1)$.
This means that every f\/inite-dimensional irreducible $\mathcal{B}$-module of type~$(1,1)$ comes from some tensor product
of evaluation modules for~$U_q(L(\mathfrak{sl}_2))$.

Let~$V$ be a~f\/inite-dimensional irreducible $U_q(L(\mathfrak{sl}_2))$-module of type $(1,1)$.
Then~$V$ is irreducible as a~$\mathcal{B}$-module and so there exists an irreducible $U_q(L(\mathfrak{sl}_2))$-module
$V'=V(\ell_1,a_1)\otimes\cdots\otimes V(\ell_n,a_n)$ such that $V$, $V'$ are isomorphic as $\mathcal{B}$-modules.
By Lemma~\ref{Lemma2}, $V$, $V'$ are isomorphic as $U_q(L(\mathfrak{sl}_2))$-modules.
This completes the proof of part~(iii) of the theorem for $U_q(L(\mathfrak{sl}_2))$.

The proof of Theorem~\ref{Theorem1} can be given by an argument very similar to the one we have seen above for the case of
$U_q(L(\mathfrak{sl}_2))$.
We prepare two more lemmas for the case of $U'_q(L(\mathfrak{sl}_2))$ to make the point clearer.
Set $(\varepsilon,\varepsilon^*)=(1,0)$ and let $\mathcal{T}_q$ be the augmented TD-algebra def\/ined by~(TD)$_0$,
(TD)$_1$ in~\eqref{Eq2},~\eqref{Eq3}.
For $s\in\mathbb{C}^\times$, let $\varphi_s$ be the embedding of $\mathcal{T}_q$ into $U'_q(L(\mathfrak{sl}_2))$ given
by~\eqref{Eq5},~\eqref{Eq6}.

\begin{lemma}\label{Lemma3}
Let $V_1$, $V_2$ be finite-dimensional irreducible $U'_q(L(\mathfrak{sl}_2))$-modules.
If $V_1$, $V_2$ are isomorphic as $\varphi_s(\mathcal{T}_q)$-modules for some $s\in\mathbb{C}^\times$, then $V_1$, $V_2$ are
isomorphic as $U'_q(L(\mathfrak{sl}_2))$-modules.
\end{lemma}

\begin{proof}
By~\eqref{Eq5}, $\varphi_s(\mathcal{T}_q)$ is generated by $se_0^++s^{-1}e_1^-k_1$, $e_1^+$ and $k_i^{\pm1}$    $(i=0,1)$.
Since $\langle e_1^\pm$, $k_1^{\pm1}\rangle$ is isomorphic to the quantum algebra $U_q(\mathfrak{sl}_2)$, the action of
$e_1^-$ on $V_i$, $i=1,2$, is uniquely determined by those of $e_1^+$, $k_1^{\pm1}\in\varphi_s(\mathcal{T}_q)$ by Lemma~\ref{Lemma2}.
Apparently the action of $e_0^+$ on $V_i$, $i=1,2$, is uniquely determined by those of $se_0^++s^{-1}e_1^-k_1$, $e_1^-$, $k_1$, and hence by that of
$\varphi_s(\mathcal{T}_q)$.
So the action of $U'_q(L(\mathfrak{sl}_2))$ on $V_i$, $i=1,2$, is uniquely determined by that of $\varphi_s(\mathcal{T}_q)$.
\end{proof}

\begin{lemma}\label{Lemma4}
Let~$V$ be a~finite-dimensional irreducible $U'_q(L(\mathfrak{sl}_2))$-module of type $(1,1)$.
Then there exists a~finite set~$\Lambda$ of nonzero scalars such that~$V$ is irreducible as
a~$\varphi_s(\mathcal{T}_q)$-module for each $s\in\mathbb{C}^\times\setminus \Lambda$.
\end{lemma}

\begin{proof}
For $s\in\mathbb{C}^\times$, regard~$V$ as a~$\varphi_s(\mathcal{T}_q)$-module.
Let~$W$ be a~minimal $\varphi_s(\mathcal{T}_q)$-submodule of~$V$.
It is enough to show that $W=V$ holds if~$s$ avoids f\/initely many scalars.
By~\eqref{Eq10} with $s_0=1$, the eigenspace decomposition of $k_1=k_0^{-1}$ on~$V$ is
$V=\bigoplus\limits_{i=0}^{d}U_{d-i}$, $k_1|_{U_{d-i}}=q^{2i-d}$, $0\leq i\leq d$.
The subalgebra $\langle e_1^\pm$, $k_1^{\pm1}\rangle$ of $U'_q(L(\mathfrak{sl}_2))$ is isomorphic to the quantum algebra
$\mathcal{U}=U_q(\mathfrak{sl}_2)$ in~\eqref{Eq15} via the correspondence of $e_1^\pm$, $k_1^{\pm1}$ to $X^\pm$, $K^{\pm1}$.
The element $(e_1^+)^{d-2i}$ maps $U_{d-i}$ onto $U_{i}$ bijectively, $0\leq i\leq[d/2]$.
Also $(e_1^-k_1)^{d-2i}$ maps $U_{i}$ onto $U_{d-i}$ bijectively, $0\leq i\leq[d/2]$.

The element $(e_1^+)^{d-2i}$ belongs to $\varphi_s(\mathcal{T}_q)$.
So $(e_1^+)^{d-2i}W\subseteq W$.
Since the mapping $(e_1^+)^{d-2i}$: $U_{d-i}\rightarrow U_{i}$ is a~bijection,
we have $\dim(W\cap U_{d-i})\leq \dim(W\cap U_i)$,
$0\leq i\leq[d/2]$.

The element $(e_1^-k_1)^{d-2i}$ does not belong to $\varphi_s(\mathcal{T}_q)$, but $(e_0^++s^{-2}e_1^-k_1)^{d-2i}$ does.
By~\eqref{Eq9}, $(e_0^++s^{-2}e_1^-k_1)^{d-2i}$ maps $U_{i}$ to $U_{d-i}$, $0\leq i\leq[d/2]$.
We want to show it is a~bijection if~$s$ avoids f\/initely many scalars.
Identify $U_{d-i}$ with $U_{i}$ as vector spaces by the bijection $(e_1^+)^{d-2i}$ between them.
Then it makes sense to consider the determinant of a~linear map from $U_{i}$ to $U_{d-i}$.
Set $t=s^{-2}$ and expand $(e_0^++te_1^-k_1)^{d-2i}$ as
\begin{gather*}
t^{d-2i}(e_1^-k_1)^{d-2i}+\text{lower terms in}~t.
\end{gather*}
Each term of the expansion gives a~linear map from $U_{i}$ to $U_{d-i}$.
So the determinant of $(e_0^++te_1^-k_1)^{d-2i}|_{U_i}$   equals
\begin{gather}
t^{(d-2i)\dim U_i}\det (e_1^-k_1)^{d-2i}\big|_{U_i}+\text{lower terms in}~t,
\label{Eq16}
\end{gather}
and this is a~polynomial in~$t$ of degree $(d-2i)\dim U_i$, since $\textrm{det}(e_1^-k_1)^{d-2i}\big|_{U_i}\neq 0$.
Let $\Lambda_i$ be the set of nonzero~$s$ such that $t=s^{-2}$ is a~root of the polynomial in~\eqref{Eq16}.
Then if $s\in\mathbb{C}^\times\setminus\Lambda_i$, $(e_0^++s^{-2}e_1^-k_1)^{d-2i}$ maps $U_{i}$ to $U_{d-i}$ bijectively.

Set $\Lambda=\cup_{i=0}^{[d/2]}\Lambda_i$.
Choose $s\in\mathbb{C}^\times\setminus\Lambda$.
Then the mapping $(e_0^++s^{-2}e_1^-k_1)^{d-2i}:U_{i}\to U_{d-i}$ is a~bijection for $0\leq i\leq[d/2]$.
Since $e_0^++s^{-2}e_1^-k_1$ belongs to $\varphi_s(\mathcal{T}_q)$, we have $(e_0^++s^{-2}e_1^-k_1)^{d-2i}W\subseteq W$
and so $\dim(W\cap U_i)\leq \dim(W\cap U_{d-i})$.
Since we have already shown $\dim(W\cap U_{d-i})\leq \dim(W\cap U_i)$, we obtain
$\dim(W\cap U_i)=\dim(W\cap U_{d-i})$,
$0\leq i\leq[d/2]$.
Therefore by Lemma~\ref{Lemma1}, we have $e_1^-W\subseteq W$.
Since $(e_0^++s^{-2}e_1^-k_1)W\subseteq W$, the inclusion $e_0^+W\subseteq W$ follows from $e_1^-W\subseteq W$ and
so~$W$ is $U'_q(L(\mathfrak{sl}_2))$-invariant.
Thus $W=V$ holds by the irreducibility of~$V$ as a~$U'_q(L(\mathfrak{sl}_2))$-module.
\end{proof}
\begin{proof}[Proof of Theorem~\ref{Theorem1}.]
We use the classif\/ication of f\/inite-dimensional irreducible $\mathcal{T}_q$-modules in the case of $(\varepsilon,
\varepsilon^*)=(1,0)$ \cite[Theorem~1.18]{augTD}:
\begin{enumerate}\itemsep=0pt
\item[(i)] A~tensor product $V=V(\ell)\otimes V(\ell_1,a_1)\otimes\dots\otimes V(\ell_n,a_n)$ of evaluation modules is
irreducible as a~$\varphi_s(\mathcal{T}_q)$-module if and only if $-s^{-2}\notin S(\ell_i,a_i)$ for all $i\in\{1,\ldots,n\}$
and $S(\ell_i,a_i)$, $S(\ell_j,a_j)$ are in general position for all $i,j\in\{1,\ldots,n\}$.
\item[(ii)] Consider two tensor products $V=V(\ell)\otimes V(\ell_1,a_1)\otimes\dots\otimes V(\ell_n,a_n)$,
$V'=V(\ell')\otimes V(\ell'_1,a'_1)\otimes\dots\otimes V(\ell'_{m},a'_{m})$ of evaluation modules and assume that they
are both irreducible as a~$\varphi_s(\mathcal{T}_q)$-module.
Then $V$, $V'$ are isomorphic as $\varphi_s(\mathcal{T}_q)$-modules if and only if $\ell=\ell'$, $n=m$ and
$(\ell_i,a_i)=(\ell'_i,a'_i)$ for all $i\in\{1,\ldots,n\}$ with a~suitable reordering of the evaluation modules
$V(\ell_1,a_1),\ldots,V(\ell_n,a_n)$.
\item[(iii)] Every f\/inite-dimensional irreducible $\mathcal{T}_q$-module $V$, $\dim V\geq 2$,
is isomorphic to a~$\mathcal{T}_q$-module $V'=V(\ell)\otimes V(\ell_1,a_1)\otimes\cdots\otimes
V(\ell_n,a_n)$ on which $\mathcal{T}_q$ acts via some embedding $\varphi_s:\mathcal{T}_q\to
U'_q(L(\mathfrak{sl}_2))$.
\end{enumerate}
Part~(i) of Theorem~\ref{Theorem1} follows immediately from the part~(i) above, due to Lemma~\ref{Lemma4}.
Part~(ii) of Theorem~\ref{Theorem1} follows immediately from the part~(ii) above, the `if' part due to Lemma~\ref{Lemma3}
(and Lemma~\ref{Lemma4})
and the `only if' part due to Lemma~\ref{Lemma4}.

We want to show part~(iii) of Theorem~\ref{Theorem1}.
Let~$V$ be a~f\/inite-dimensional irreducible $U'_q(L(\mathfrak{sl}_2))$-module of type $(1,1)$.
By Lemma~\ref{Lemma4}, there exists a~nonzero scalar~$s$ such that~$V$ is irreducible as a~$\varphi_s(\mathcal{T}_q)$-module.
By the part~(iii) above, for the proof of which Drinfel'd polynomials play the key role, the $\mathcal{T}_q$-module~$V$
via $\varphi_s$ is isomorphic to some $\mathcal{T}_q$-module $V'=V(\ell)\otimes V(\ell_1,a_1)\otimes\dots\otimes
V(\ell_n,a_n)$ via some embedding $\varphi_{s'}$ of $\mathcal{T}_q$ into $U'_q(L(\mathfrak{sl}_2))$.
Since $k_0$ has the same eigenvalues on $V$, $V'$, we have $s=s'$ and so $V$, $V'$ are isomorphic as
$\varphi_s(\mathcal{T}_q)$-modules.
By Lemma~\ref{Lemma3}, $V$, $V'$ are isomorphic as $U'_q(L(\mathfrak{sl}_2))$-modules.
Part~(iv) will be shown in the next section.
\end{proof}

\section{Intertwiners: Proof of Theorem~\ref{Theorem1}(iv)}\label{Section4}

In this section, we show that for $\ell$, $m\in\mathbb{Z}_{>0}$, $a\in\mathbb{C}^\times$, there exists an intertwiner between
the $U'_q(L(\mathfrak{sl}_2))$-modules $V(\ell,a)\otimes V(m)$, $V(m)\otimes V(\ell,a)$, i.e., a~nonzero linear map~$R$
from $V(\ell,a)\otimes V(m)$ to $V(m)\otimes V(\ell,a)$ such that
\begin{gather}
R\Delta(\xi)=\Delta(\xi)R,
\qquad
\forall\,\xi\in U'_q(L(\mathfrak{sl}_2)),
\label{Eq17}
\end{gather}
where~$\Delta$ is the coproduct from~\eqref{Eq12}.
If such an intertwiner~$R$ exists, then it is routinely concluded that $V(\ell,a)\otimes V(m)$ is isomorphic to
$V(m)\otimes V(\ell,a)$ as $U'_q(L(\mathfrak{sl}_2))$-modules and any other intertwiner is a~scalar multiple of~$R$,
since $V(m)\otimes V(\ell,a)$ is irreducible as a~$U'_q(L(\mathfrak{sl}_2))$-module by Theorem~\ref{Theorem1}.

Using the theory of Drinfel'd polynomials~\cite{augTD}
for the augmented TD-algebra $\mathcal{T}_q=\mathcal{T}_q^{(\varepsilon,
\varepsilon^*)}$ with $(\varepsilon, \varepsilon^*)=(1,0)$, we shall f\/irstly show that $V(\ell,a)\otimes V(m)$ is
isomorphic to $V(m)\otimes V(\ell,a)$ as $U'_q(L(\mathfrak{sl}_2))$-modules.
This proves Theorem~\ref{Theorem1}(iv), since it is well-known~\cite[Theorem~4.11]{CP} that $V(\ell_i, a_i)\otimes V(\ell_j,a_j)$
and $V(\ell_j, a_j)\otimes V(\ell_i,a_i)$ are isomorphic as $U_q(\widehat{\mathfrak{sl}}_2)$-modules, if $S(\ell_i, a_i)$
and $S(\ell_j, a_j)$ are in general position.
Finally we shall construct an intertwiner explicitly.

Let us denote the $U'_q(L(\mathfrak{sl}_2))$-modules $V(\ell,a)\otimes V(m)$, $V(m)\otimes V(\ell,a)$ by $V$, $V'$:
\begin{gather*}
V=V(\ell,a)\otimes V(m),
\qquad
V'=V(m)\otimes V(\ell,a).
\end{gather*}
Recall we assume that the integers $\ell$, $m$ and the scalar~$a$ are nonzero.
Let us denote a~standard basis of the $U'_q(L(\mathfrak{sl}_2))$-module $V(\ell,a)$
(resp.~$V(m)$) by $v_0, v_1, \ldots, v_\ell$ (resp.~$v'_0, v'_1, \ldots, v'_m$) in the sense of~\eqref{Eq11}.
Recall~$V(m)$ is an abbreviation of $V(m,0)$ and the action of $U'_q(L(\mathfrak{sl}_2))$ on $V$, $V'$ are via the
coproduct~$\Delta$ of~\eqref{Eq12}.

Let $\mathcal{U}$ denote the subalgebra of $U'_q(L(\mathfrak{sl}_2))$ generated by $e_1^\pm$, $k_1^\pm$.
The subalgebra $\mathcal{U}$ is isomorphic to the quantum algebra $U_q(\mathfrak{sl}_2)$.
Let $V(n)$ denote the irreducible $\mathcal{U}$-module of dimension $n+1$: $V(n)$ has a~standard basis
$x_0,x_1,\ldots,x_n$ on which $\mathcal{U}$ acts as
\begin{gather*}
 k_1x_i=q^{n-2i}x_i,
\qquad
e^+_1x_i=[n-i+1]x_{i-1},
\qquad
e^-_1x_i=[i+1]x_{i+1},
\end{gather*}
where $[t]=[t]_q=(q^t-q^{-t})/(q-q^{-1})$, $x_{-1}=x_{n+1}=0$.
We call~$x_n$ (resp.~$x_0$) the \textit{lowest} (\textit{highest}) \textit{weight vector}:
$k_1x_n=q^{-n}x_n$, $e^-_1x_n=0$ $(k_1x_0=q^nx_0$, $e^+_1x_0=0)$.
Note that $V(\ell,a)$ is isomorphic to $V(\ell)$ as $\mathcal{U}$-modules.

By the Clebsch--Gordan formula, $V=V(\ell,a)\otimes V(m)$ is decomposed into the direct sum of $\mathcal{U}$-submodules
$\widetilde{V}(n)$, $|\ell-m|\leq n\leq \ell+m$, $n\equiv\ell+m\bmod 2$,
where $\widetilde{V}(n)$ is the unique irreducible
$\mathcal{U}$-submodule of~$V$ isomorphic to $V(n)$.
With $n=\ell+m-2\nu$, we have
\begin{gather}
V=V(\ell,a)\otimes V(m)=\bigoplus_{\nu=0}^{\min \{\ell,m\}}\widetilde{V}(\ell+m-2\nu).
\label{Eq19}
\end{gather}
Let $\widetilde{x}_n$ denote a~lowest weight vector of the $\mathcal{U}$-module $\widetilde{V}(n)$.
So
\begin{gather*}
\Delta(k_1)\widetilde{x}_n=q^{-n}\widetilde{x}_n,
\qquad
\Delta(e_1^-)\widetilde{x}_n=0.
\end{gather*}
Since~$V$ has a~basis $\{v_{\ell-i}\otimes v'_{m-j}\,|\,  0\leq i\leq\ell,0\leq j\leq m\}$ and $k_1$ acts on it~by
$\Delta(k_1)(v_{\ell-i}\otimes v'_{m-j})=q^{-(\ell+m)+2(i+j)}v_{\ell-i}\otimes v'_{m-j}$, the lowest weight vector
$\widetilde{x}_n$ of $\widetilde{V}(n)$ can be expressed as
\begin{gather*}
\widetilde{x}_n=\sum\limits_{i+j=\nu}c_jv_{\ell-i}\otimes v'_{m-j},
\qquad
n=\ell+m-2\nu.
\end{gather*}
Solving $\Delta(e_1^-)\widetilde{x}_n=0$ for the coef\/f\/icients $c_j$, we obtain
\begin{gather*}
\frac{c_j}{c_{j-1}}=-q^{m-2j+2}\frac{[\ell-\nu+j]}{[m-j+1]}
\end{gather*}
and so with a~suitable choice of $c_0$
\begin{gather}
\widetilde{x}_n=\sum\limits_{j=0}^\nu(-1)^jq^{j(m-j+1)}[\ell-\nu+j]![m-j]!v_{\ell-\nu+j}\otimes v'_{m-j},
\label{Eq22}
\end{gather}
where $n=\ell+m-2\nu$ and $[t]!=[t][t-1]\cdots[1]$.

\begin{lemma}\label{Lemma5}
$\Delta(e_0^+)\widetilde{x}_n=aq\widetilde{x}_{n+2}$.
\end{lemma}

\begin{proof}
By~\eqref{Eq12}, we have $\Delta(e_0^+)=e_0^+\otimes 1+k_0\otimes e_0^+$.
By~\eqref{Eq11}, the element $e_0^+$ vanishes on $V(m)$ and acts on $V(\ell,a)$ as $aqe_1^-$.
Since $e^-_1v_{\ell-\nu+j}=[\ell-(\nu -1)+j]v_{\ell-(\nu -1)+j}$, the result follows from~\eqref{Eq22}, using $v_{\ell+1}=0$.
\end{proof}

\begin{corollary}\label{Cor1}
Any nonzero $U'_q(L(\mathfrak{sl}_2))$-submodule of $V(\ell,a)\otimes V(m)$ contains $\widetilde{x}_{\ell+m}$, the
lowest weight vector of the $\mathcal{U}$-module $V(\ell,a)\otimes V(m)$.
\end{corollary}

We are ready to prove the following

\begin{theorem}\label{Theorem2}
The $U'_q(L(\mathfrak{sl}_2))$-modules $V(\ell,a)\otimes V(m)$, $V(m)\otimes V(\ell,a)$ are isomorphic for every
$\ell,m\in\mathbb{Z}_{>0}$, $a\in\mathbb{C}^\times$.
\end{theorem}

\begin{proof}
Let $\mathcal{T}_q=\mathcal{T}_q^{(\varepsilon, \varepsilon^*)}$ be the augmented TD-algebra with $(\varepsilon,
\varepsilon^*)=(1,0)$.
Let $\varphi_s: \mathcal{T}_q\to U'_q(L(\mathfrak{sl}_2))$ denote the embedding of $\mathcal{T}_q$ into
$U'_q(L(\mathfrak{sl}_2))$ given in~\eqref{Eq6}.
By \cite[Theorem~5.2]{augTD}, the Drinfel'd polynomial $P_V(\lambda)$ of the $\varphi_s(\mathcal{T}_q)$-module
$V=V(\ell,a)\otimes V(m)$ turns out to be
\begin{gather*}
P_V(\lambda)=\lambda^m\prod\limits_{i=0}^{\ell-1}\big(\lambda+aq^{2i-\ell+1}\big).
\end{gather*}
(Note that the parameter~$s$ of the embedding $\varphi_s$ does not appear in $P_V(\lambda)$.
So the polynomial~$P_V(\lambda)$ can be called the Drinfel'd polynomial attached to the
$U'_q(L(\mathfrak{sl}_2))$-module~$V$.)

Let~$W$ be a~minimal~$U'_q(L(\mathfrak{sl}_2))$-submodule of $V=V(\ell,a)\otimes V(m)$; notice that we have shown the
irreducibility of the $U'_q(L(\mathfrak{sl}_2))$-module $V'=V(m)\otimes V(\ell,a)$ in Theorem~\ref{Theorem1} but not yet of
$V=V(\ell,a)\otimes V(m)$.
By Corollary~\ref{Cor1}, $W$ contains the lowest and hence highest weight vectors of~$V$.
In particular, the irreducible $U'_q(L(\mathfrak{sl}_2))$-module~$W$ is of type $(1, 1)$.
By Lemma~\ref{Lemma4}, there exists a~f\/inite set~$\Lambda$ of nonzero scalars such that~$W$ is irreducible as
a~$\varphi_s(\mathcal{T}_q)$-module for any $s\in\mathbb{C}^\times\setminus\Lambda$.
By the def\/inition, the Drinfel'd polynomial~$P_W(\lambda)$ of the irreducible $\varphi_s(\mathcal{T}_q)$-module~$W$
coincides with $P_V(\lambda)$:
\begin{gather*}
P_W(\lambda)=P_V(\lambda).
\end{gather*}

By Theorem~\ref{Theorem1}, $V'=V(m)\otimes V(\ell,a)$ is irreducible as a~$U'_q(L(\mathfrak{sl}_2))$-module.
So by Lemma~\ref{Lemma4}, there exists a~f\/inite set $\Lambda'$ of nonzero scalars such that $V'$ is irreducible as
a~$\varphi_s(\mathcal{T}_q)$-module for any $s\in\mathbb{C}^\times\setminus \Lambda'$.
By \cite[Theorem~5.2]{augTD}, the Drinfel'd polynomial~$P_{V'}(\lambda)$ of the irreducible $\varphi_s(\mathcal{T}_q)$-module
$V'$ coincides with $P_V(\lambda)$:
\begin{gather*}
P_{V'}(\lambda)=P_V(\lambda).
\end{gather*}

Both of the irreducible $\varphi_s(\mathcal{T}_q)$-modules $W$, $V'$ have type~$s$, diameter $d=\ell+m$ and the Drinfel'd
polynomial~$P_V(\lambda)$.
By \cite[Theorem 1.9$'$]{augTD}, $W$ and $V'$ are isomorphic as $\varphi_s(\mathcal{T}_q)$-modules.
By Lemma~\ref{Lemma3}, $W$ and $V'$ are isomorphic as $U'_q(L(\mathfrak{sl}_2))$-modules.
In particular, $\dim W=\dim V'$.
Since $\dim V'=\dim V$, we have $W=V$, i.e., $V$ and $V'$ are isomorphic as $U'_q(L(\mathfrak{sl}_2))$-modules.
\end{proof}

Finally we want to construct an intertwiner~$R$ between the irreducible $U'_q(L(\mathfrak{sl}_2))$-mo\-du\-les~$V$,~$V'$.
Regard $V'=V(m)\otimes V(\ell,a)$ as a~$\mathcal{U}$-module.
By the Clebsch--Gordan formula, we have the direct sum decomposition
\begin{gather}
V'=V(m)\otimes V(\ell,a)=\bigoplus_{\nu=0}^{\min \{\ell,m\}}\widetilde{V}'(\ell+m-2\nu),
\label{Eq26}
\end{gather}
where $\widetilde{V}'(n)$ is the unique irreducible $\mathcal{U}$-submodule of $V'$ isomorphic to $V(n)$, $n=\ell+m-2\nu$.
Let $\widetilde{x}'_n$ be a~lowest weight vector of the $\mathcal{U}$-module $\widetilde{V}'(n)$.
By~\eqref{Eq22}, we have
\begin{gather}
\widetilde{x}'_n=\sum\limits_{j=0}^\nu(-1)^jq^{j(\ell-j+1)}[m-\nu+j]![\ell-j]!v'_{m-\nu+j}\otimes v_{\ell-j}
\label{Eq27}
\end{gather}
up to a~scalar multiple, where $n=\ell+m-2\nu$.
It can be easily checked as in Lemma~\ref{Lemma5} that the lowest weight vectors $\widetilde{x}'_n$, $n=\ell+m-2\nu$,
$0\leq\nu\leq\min \{\ell,m\}$,
are related~by
\begin{gather}
(e_1^-\otimes 1)\widetilde{x}'_n=\widetilde{x}'_{n+2},
\label{Eq28}
\end{gather}
where $V'=V(m)\otimes V(\ell,a)$ is regarded as a~$(\mathcal{U}\otimes\mathcal{U})$-module in the natural way.

\begin{lemma}\label{Lemma6}
$\Delta(e_0^+)\widetilde{x}'_n=-aq\cdot q^{n+2}\widetilde{x}'_{n+2}$.
\end{lemma}

\begin{proof}
We have $\Delta(e_0^+)\widetilde{x}'_n=aq(k_1^{-1}\otimes e_1^-)\widetilde{x}'_n$, since
$\Delta(e_0^+)=e_0^+\otimes 1+k_0\otimes e_0^+$, and $e_0^+$ vanishes on $V(m)$ and acts on $V(\ell,a)$ as
$aqe_1^-$.
Express $k_1^{-1}\otimes e_1^-$ as $k_1^{-1}\otimes e_1^- =(k_1^{-1}\otimes 1)(1\otimes e_1^-)
=(k_1^{-1}\otimes 1)(\Delta(e_1^-)-e_1^-\otimes k_1^{-1}) =(k_1^{-1}\otimes 1)\Delta(e_1^-)-k_1^{-1}e_1^-\otimes k_1^{-1}
=(k_1^{-1}\otimes 1)\Delta(e_1^-)-q^2(e_1^-\otimes 1)\Delta(k_1^{-1})$. Since
$\Delta(e_1^-)\widetilde{x}'_n=0$, $\Delta(k_1^{-1})\widetilde{x}'_n=q^n\widetilde{x}'_n$, the result follows from~\eqref{Eq28}.
\end{proof}

There exists a~unique linear map
\begin{gather*}
R_n: \ V=V(\ell,a)\otimes V(m)\to\widetilde{V}'(n)
\end{gather*}
that commutes with the action of $\mathcal{U}$ and sends $\widetilde{x}_n$ to $\widetilde{x}'_n$.
The linear map $R_n$ vanishes on $\widetilde{V}(t)$ for $t \neq n$ and af\/fords an isomorphism between $\widetilde{V}(n)$
and $\widetilde{V}'(n)$ as $\mathcal{U}$-modules.
If~$R$ is an intertwiner in the sense of~\eqref{Eq17}, then
$R$ can be expressed as
\begin{gather}
R=\sum\limits_{\nu=0}^{\min \{\ell,m\}}\alpha_\nu R_{\ell+m-2\nu},
\label{Eq29}
\end{gather}
regarding~$R$ as an intertwiner for the $\mathcal{U}$-modules $V$, $V'$.
By~\eqref{Eq17}, we have
\begin{gather}
R\Delta(e_0^+)=\Delta(e_0^+)R.
\label{Eq30}
\end{gather}
Apply~\eqref{Eq30} to the lowest weight vector $\widetilde{x}_n$ in~\eqref{Eq22}.
By Lemma~\ref{Lemma5}, $\Delta(e_0^+)\widetilde{x}_n=aq\widetilde{x}_{n+2}$ and so with $n=\ell+m-2\nu$, we have
\begin{gather}
R\Delta(e_0^+)\widetilde{x}_n=aq\alpha_{\nu-1}\widetilde{x}'_{n+2}.
\label{Eq31}
\end{gather}
On the other hand, $R\widetilde{x}_n=\alpha_\nu\widetilde{x}'_n$, $n=\ell+m-2\nu$, and so by Lemma~\ref{Lemma6}, we have
\begin{gather}
\Delta(e_0^+)R\widetilde{x}_n=-aq\alpha_\nu q^{n+2}\widetilde{x}'_{n+2}.
\label{Eq32}
\end{gather}
By~\eqref{Eq31},~\eqref{Eq32}, we have $\alpha_\nu/\alpha_{\nu-1}=-q^{-n-2}=-q^{-\ell-m+2(\nu-1)}$ and so
\begin{gather}
\alpha_\nu=(-1)^\nu q^{-\nu(\ell+m-\nu+1)}
\label{Eq33}
\end{gather}
by choosing $\alpha_0=1$.
An intertwiner exists by Theorem~\ref{Theorem2}.
If it exists, it has to be in the form of~\eqref{Eq29},~\eqref{Eq33}.
Thus we obtain the following.
\begin{theorem}\label{Theorem3}
The linear map
\begin{gather*}
R=\sum\limits_{\nu=0}^{\min \{\ell,m\}}(-1)^\nu q^{-\nu(\ell+m-\nu+1)}R_{\ell+m-2\nu}
\end{gather*}
is an intertwiner between the $U'_q(L(\mathfrak{sl}_2))$-modules $V(\ell,a)\otimes V(m)$, $V(m)\otimes V(\ell,a)$.
Any other intertwiner is a~scalar multiple of~$R$.
\end{theorem}

\begin{remark}
Let $R(a,b)$ be an intertwiner between the irreducible $U'_q(L(\mathfrak{sl}_2))$-modules $V=V(\ell,a)\otimes
V(m,b)$, $V'=V(m,b)\otimes V(\ell,a)$, where $a \neq 0$, $b \neq 0$:
\begin{gather*}
R(a,b): \ V=V(\ell,a)\otimes V(m,b)\to V'=V(m,b)\otimes V(\ell,a).
\end{gather*}
As in~\eqref{Eq29}, we write
\begin{gather*}
R(a,b)=\sum\limits_{\nu=0}^{\min \{\ell,m\}}\alpha_\nu R_{\ell+m-2\nu}.
\end{gather*}
Recall~$R_n$ is the linear map from~$V$ to $V'$ that commutes with the action of $\mathcal{U}=\langle e_1^\pm$, $k_1^\pm
\rangle$ and sends $\widetilde{x}_n$ to $\widetilde{x}'_n$, where $\widetilde{x}_n$, $\widetilde{x}'_n$ are the lowest
weight vectors of $\widetilde{V}(n)$, $\widetilde{V}'(n)$ from~\eqref{Eq19},~\eqref{Eq26} that are explicitly given by~\eqref{Eq22},~\eqref{Eq27}
and satisfy $(e_1^-\otimes1)\widetilde{x}_n=\widetilde{x}_{n+2}$, $(e_1^-\otimes1)\widetilde{x}'_n=\widetilde{x}'_{n+2}$
as in~\eqref{Eq28}.
Since the $U'_q(L(\mathfrak{sl}_2))$-modules $V$, $V'$ can be extended to $U_q(\widehat{\mathfrak{sl}}_2)$-modules, we
have by~\cite[Theorem~5.4]{CP}
\begin{gather}
\alpha_{\nu} = \prod\limits_{j=0}^{\nu -1}\frac{a-bq^{\ell+m-2j}}{b-aq^{\ell+m-2j}},
\label{Eq34}
\end{gather}
where we choose $\alpha_0=1$.
Note that the denominator and the numerator of~\eqref{Eq34} are non-zero, since $V$, $V'$ are assumed to be irreducible and so
$S(\ell,a)$, $S(m,b)$ are in general position.

The intertwiner $R(a,b)$ with $a\neq 0$, $b \neq 0$ is derived from the universal $R$-matrix for the quantum af\/f\/ine
algebra $U_q(\widehat{\mathfrak{sl}}_2)$~\cite{TK}.
If we put $b=0$ in~\eqref{Eq34}, then the spectral parameter~$u$ disappears, where $q^{2u}=a/b$, and we get~\eqref{Eq33}.
In this sense, the intertwiner~$R$ of Theorem~\ref{Theorem3} is related to the universal $R$-matrix for
$U_q(\widehat{\mathfrak{sl}}_2)$, but we cannot expect that~$R$ comes from it directly, because the
$U'_q(L(\mathfrak{sl}_2))$-modules $V=V(\ell,a)\otimes V(m)$, $V'=V(m)\otimes V(\ell,a)$ cannot be extended to
$U_q(\widehat{\mathfrak{sl}}_2)$-modules.
In order to derive both of the intertwiners $R(a,b)$, $R$ from a~universal $R$-matrix directly, we need to construct it
for the subalgebra $U_q'(\widehat{\mathfrak{sl}}_2)$ of $U_q(\widehat{\mathfrak{sl}}_2)$.
\end{remark}

\subsection*{Acknowledgements}

The research of the second author is supported by Science Foundation of Anhui University (Grant No.~J10117700037).

\pdfbookmark[1]{References}{ref}
\LastPageEnding

\end{document}